\documentclass[12pt]{article}
\usepackage[final]{epsfig}
\usepackage{graphics}
\usepackage{amsmath}
\usepackage{amsfonts}
\usepackage{latexsym}
\usepackage{amssymb}
\usepackage{graphicx}
\usepackage{epstopdf}
\newtheorem{lemma}{Lemma}
\newtheorem{proposition}[lemma]{Proposition}
\newtheorem{remark}[lemma]{Remark}

\newtheorem{theorem}[lemma]{Theorem}

\begin{document}
\newcommand{\eps}{{\varepsilon}}
\newcommand{\proofend}{$\Box$\bigskip}
\newcommand{\C}{{\mathbb C}}
\newcommand{\Q}{{\mathbb Q}}
\newcommand{\R}{{\mathbb R}}
\newcommand{\Z}{{\mathbb Z}}
\newcommand{\RP}{{\mathbb{RP}}}
\newcommand{\CP}{{\mathbb{CP}}}

\def\id{\mathop{\rm id}}
\def\Ker{\mathop{\rm Ker}}
\def\rank{\mathop{\rm rank}}

\title {Self-dual polygons and self-dual curves}
\author{Dmitry Fuchs\thanks{
Department of Mathematics,
University of California, Davis, CA 95616, USA;
e-mail: \tt{fuchs@math.ucdavis.edu}
}
\ and Serge Tabachnikov\thanks{
Department of Mathematics,
Pennsylvania State University, University Park, PA 16802, USA;
e-mail: \tt{tabachni@math.psu.edu}
}
\\
}
\date{\today}
\maketitle

\section{Introduction} \label{intro}

The projective plane $P$ is the projectivization of 3-dimensional space $V$ (we consider the cases of two ground fields, $\R$ and $\C$); the points  and the lines of $P$ are 1- and 2-dimensional subspaces of $V$. The dual projective plane $P^\ast$ is the projectivization of the dual space $V^\ast$. Assign to a subspace in $V$ its annihilator in $V^\ast$. This gives a correspondence between the points in $P$ and the lines in $P^\ast$, and between the lines in $P$ and the points in $P^\ast$, called the projective duality. Projective duality preserves the incidence relation:
if a point $A$ belongs to a line $B$ in $P$ then the dual point $B^\ast$ belongs to the dual line $A^\ast$ in $P^\ast$.  Projective duality is an involution: $(A^\ast)^\ast=A$.

Projective duality extends to polygons in $P$. An $n$-gon is a cyclically ordered collection of $n$ points and $n$ lines satisfying the incidences: two consecutive vertices lie on the respective side, and two consecutive sides pass through the respective vertex. We assume that our polygons are non-degenerate: no three consecutive vertices are collinear. Thus to every polygon $L\subset P$ there corresponds the dual polygon $L^\ast \subset P^\ast$. A polygon $L$ is called self-dual if there exists a projective map $P \to P^\ast$ that takes $L$ to $L^\ast$.

Projective duality also extends to locally convex smooth curves. A smooth curve $\gamma\subset P$ determines a one-parameter family of its tangent lines, and projective duality takes it to a one-parameter family of points in $P^\ast$, the dual curve $\gamma^\ast \subset P^\ast$. If $\gamma$ is locally convex then $\gamma^\ast$ is smooth as well. One has $(\gamma^\ast)^\ast=\gamma$. Projective duality further extends to a broader class of curves with inflections and cusps, called {\it wave fronts} (see end of Section \ref{polycurves} for a precise definition). Projective duality interchanges inflections and cusps. One defines self-dual curves similarly to self-dual polygons. 

A motivation for this work is the following problem (No 1994-17 in \cite{Arn}) of V. Arnold:

\begin{quote}
Find all projective curves equivalent to their duals. {\it The answer seems to be unknown even 
in $\RP^2$}.
\end{quote}

(A traditional interpretation of this question would be to consider  algebraic curves, in which case the Plucker formulas play a critical role; see \cite{Hol}. In particular, the Plucker formulas imply that a {\it non-singular} self-dual algebraic curve is a conic; however, there exist other, singular, self-dual curves, for example, $y=x^3$, projectively equivalent to its dual $y=x^{3/2}$.)

The main result of this paper is a description of self-dual polygons in $\CP^2$. Let $A_1,A_3,\dots,A_{2n-1}\in P$ (where the indices are odd residues modulo $2n$) be the vertices of an $n$-gon, and let $B_2, B_4,
\dots,B_{2n}$ (where the indices are  even residues modulo $2n$) be its respective sides: $B_{2i}=A_{2i-1}A_{2i+1}$ for all $i$. Let $m$ be an odd number, $1\le m\le n$. The $n$-gon $L$ is called {\it $m$-self-dual} if there exists a projective map $g\colon P \to P^\ast$ such that $g(A_i)=B_{i+m}^\ast$ for all $i$. An example of an $m$-self-dual $n$-gon, for arbitrary $m$, is a regular $n$-gon.  Denote by 
${\mathcal M}_{m,n}$ the moduli space of $m$-self-dual $n$-gons. Our result is as follows.

\begin{theorem} \label{main}
If $(m,n)=1$ then ${\mathcal M}_{m,n}$ consists of one point, the class of a regular $n$-gon. If $m<n,\ (m,n)>1$ and $n\neq 2m$ then dim ${\mathcal M}_{m,n}=(m,n)-1$. Finally, dim ${\mathcal M}_{m,2m}=m-3$ and dim ${\mathcal M}_{n,n}=n-3$.
\end{theorem}

Note, for comparison, that the dimension of the moduli space of $n$-gons is $2n-8$.
The proof of Theorem \ref{main} occupies Sections \ref{m=n} and \ref{m<n}. These sections also contain explicit constructions of self-dual polygons. 

The map $g\colon P \to P^\ast$, associated with an $m$-self-dual $n$-gon, determines a linear map $V\to V^\ast$, defined up to a factor, and therefore a bilinear form $F$ on $V$. We prove that if the polygon is not a multiple of another polygon then $F$ is symmetric if and only if $m=n$, see Proposition \ref{symmetric}. We also show that if an $n$-self-dual $n$-gon is convex then the symmetric bilinear form $F$ is definite, Proposition \ref{convex}. 

We make additional observations. First, every pentagon is 5-self-dual, see Proposition \ref{pentagons} below (five is the first interesting number because all triangles are projectively equivalent, and so are all quadrilaterals). Secondly, if an $n$-gon with odd $n$ is inscribed in a conic and circumscribed about a conic then it is $n$-self-dual, see Proposition \ref{poncelet}. However the moduli space of such ``Poncelet" polygons has dimension two, which is less than $n-3$ for $n\geq 7$.

In the real case, one can also interpret a polygon as a closed polygonal curve. In Section \ref{polycurves},  we define a polygonal curve as a polygon in $\RP^2$ with two additional structures: every two consecutive vertices $A_{2i-1},A_{2i+1}$ partition the real projective line $B_{2i}$ into two segments, and one of these segments is chosen (as a side); every two consecutive sides $B_{2i}$ and $B_{2i+2}$ determine two pairs of vertical angles at the vertex $A_{2i+1}$, and one of these pairs is chosen (as an exterior angle). Polar duality naturally extends to these polygonal curves, and one can consider self-dual polygonal curves. A given $n$-gon $L$ gives rise to $2^{2n}$ polygonal curves.
We prove (Proposition \ref{choices}) that if $L$ is $m$-self-dual then, out of these $2^{2n}$ polygonal curves,  $2^{(m,n)}$ are $m$-self-dual.

Section \ref{sdcurves} concerns self-dual curves and wave fronts in the real projective plane $P$. We do not attempt to give a complete classification of such curves.  A curve $\gamma(t) \subset P,\ t\in S^1=\R/2\pi\Z$, is called self-dual if there exists a projective transformation $g\colon P\to P^\ast$ and a diffeomorphism $\varphi$ of $S^1$ such that $g(\gamma(\varphi (t))=\gamma^\ast(t)$. The diffeomorphism $\varphi$ is a continuous analog of the cyclic shift by $m$ in the definition of $m$-self-dual $n$-gons.

A number of results that we establish for polygons have analogs for curves. For example, the bilinear form $F$ is symmetric if and only if $\varphi^2=\id$, see Proposition \ref{curvesymm}. If, in addition, a self-dual curve is convex then $F$ is definite, Proposition \ref{curvepos}. 

We observe that curves of constant width $\displaystyle \frac\pi2$ on the unit sphere $S^2$ project to self-dual curves with $\varphi^2=\id$ in $\RP^2$. We construct such curves of constant width as Legendrian curves in the  manifold of contact elements of $S^2$ satisfying certain monodromy conditions. We also give a similar description to self-dual curves with the diffeomorphism $\varphi$ having higher order than 2. This description leads to  explicit formulas for self-dual curves.

Finally, we briefly discuss the Radon curves, the unit circles in two-dimensional normed spaces for which the orthogonality relation is symmetric. Radon curves have been extensively studied; they provide examples of projectively self-dual curves. 

Let us finish this introduction with a question: can a smooth convex self-dual curve, other than a conic, be an oval of an algebraic curve?

\bigskip

\paragraph{Acknowledgments.} Many thanks to J. C. Alvarez, V. Ovsienko and R. Schwartz for stimulating discussions. It is a pleasure to acknowledge the hospitality of the Research in Pairs program at the Mathematical Institute at Oberwolfach (MFO) where this work was done. The second  author was partially supported by an NSF grant DMS-0555803.

\section{Polygons and duality} \label{poldual}

We use the notation from Section \ref{intro}. Let $P={\C}P^2$ and $V={\C}^3$. Let $L$ be an $n$-gon in $P$ whose vertices are  $A_1,A_3,\dots,A_{2n-1}$ and whose sides are $B_2,B_4,
\dots,B_{2n}$. The dual $n$-gon $L^\ast$ in the dual projective plane $P^\ast$ has the vertices $B_2^\ast,B_4^\ast, \dots,B_{2n}^\ast$ and the sides $A_1^\ast,A_3^\ast,\dots,A_{2n-1}^\ast$. 

Assume that $L$ is $m$-self-dual where $m$ be an odd number, $1\le m\le n$. Then there exists a linear isomorphism $f\colon V\to V^\ast$ that takes the line $A_i \subset V$ to the line $B_{i+m}^\ast \subset V^\ast$ for all $i$. Along with $f$, we shall consider the corresponding projective isomorphism $\hat f\colon P\to P^\ast$ and the bilinear form $F$ on $V$, $F(v,w)=\langle f(v),w\rangle$. Obviously, for a given 
$m$-self-dual polygon, $\hat f$ is unique, while $f$ and $F$ are unique up to a non-zero constant factor.

Along with an $n$-gon $L=A_1A_3,\dots,A_{2n-1}$, we can consider the $kn$-gon $kL=A_1A_3\dots
A_{2kn-1}$ with $A_i=A_{i+2n}$. Obviously, $(kL)^\ast=kL^\ast$, and if $L$ is $m$-self-dual, then $kL$
is $(m+2rn)$-self-dual for $r=0,1,\dots,k-1$. A polygon $L$ is called {\it simple}, if $L\ne kL'$ for any 
$k>1$ and $L'$. 

\begin{proposition}\label{symmetric}
Let $L$ be a simple $m$-self-dual $n$-gon and $f\colon V\to V^\ast$ be the corresponding isomorphism. Then $f$ is self-adjoint $($or, equivalently, $F$ is symmetric$)$ if and only if $m=n$.
\end{proposition}

\paragraph{\it Proof.} Notice that \[\begin{array} {rl} F(A_i,A_j)=0&\Leftrightarrow \langle 
F(A_i),A_j\rangle=0\Leftrightarrow\langle B^\ast_{i+m},A_j\rangle=0\\ &\Leftrightarrow A_j\in 
B_{i+m}=A_{i+m-1}A_{i+m+1}.\end{array}\]  In particular, $F(A_i,A_{i+m\pm1})=0$ for all $i$.

Let $F$ be symmetric. Then, for all $i$,  $F(A_{i+m\pm 1},A_i)=0$, and hence
\[\begin{array} {rl} A_i\in&A_{(i+m-1)+(m-1)}A_{(i+m-1)+(m+1)}=A_{i+2m-2}A_{i+2m},\\  A_i\in&A_{(i+m+1)+(m-1)}A_{(i+m+1)+(m+1)}=A_{i+2m}A_{i+2m+2}.\end{array}\]Since the points 
$A_{i+2m-2}, A_{i+2m}, A_{i+2m+2}$ are not collinear, this means that $A_i=A_{i+2m}$. Hence $m=n$ (the polygon $L$ is simple!).

Let $m=n$.  For every $i$, $F(A_i,A_{i+m\pm1})=0$, and, in addition to that, 
\[ \begin{array} {rl} F(A_{i+m-1},A_i)&=F(A_{i+m-1},A_{i+2m})=
F(A_{i+m-1},A_{(i+m-1)+(m+1)})=0,\\ F(A_{i+m+1},A_i)&=
F(A_{i+m+1},A_{i+2m})=F(A_{i+m+1},A_{(i+m+1)+(m-1)})=0.\end{array}\]
This implies that the linear forms $F(A_i,-)$ and $F(-,A_i)$ are proportional for
every $i$, that is, there exist non-zero complex numbers $\lambda_i$ such that $ F(A_i,x)=\lambda_i F(x,A_i)$ for all $i$ and $x$. Hence $ F(A_i,A_j)=\lambda_i F(A_j,A_i)=\lambda_i\lambda_j F(A_i,A_j)$, so if $ F(A_i,A_j)\ne0$, then $\lambda_i=\lambda_j^{-1}$. But $ F(A_i,A_{i+m-3})\ne0$ (because $A_{i+m-3}$ does not belong to the line $A_{i+m-1}A_{i+m+1}$, which is the zero locus of the form $ F(A_i,x)$). Hence $$\lambda_i=\lambda_{i+m-3}^{-1}=\lambda_{i+2(m-3)}=\lambda_{i+3(m-3)}^{-1}=\dots=\lambda_{i+m(m-3)}^{-1}=
\lambda_i^{-1},$$so $\lambda_i=\pm1$. We state that all $\lambda_i$'s are the same. Indeed, if $\lambda_i=1,\lambda_j=-1$ for some $i,j$, then for every $k$, one of $ F(A_i,A_k), 
 F(A_j,A_k)$ must be 0, that is, $A_k$ belongs to one of two lines, $A_{i+m-1}A_{i+m+1}$ and 
$A_{j+m-1}A_{j+m+1}$, that is, all vertices of the polygon belong to two lines, which is impossible, since the three lines $A_1A_2,A_2A_3,A_3A_4$ are all different. We see that our form $ F$ is either symmetric or skew symmetric. However it cannot be skew-symmetric because it is non-degenerate, and all skew-symmetric forms in an odd-dimensional space are degenerate. \proofend

\begin{figure}[hbtp]
\centering
\includegraphics[width=2in]{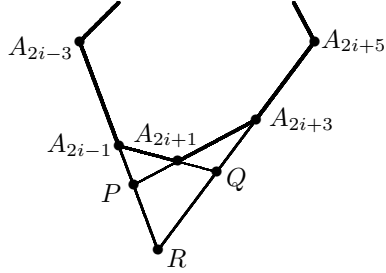}
\caption{Cross-ratios at the vertices of a polygon}
\label{cross}
\end{figure}

\begin{remark} \label{coord}
{\rm The class of projective equivalence of an $n$-gon $L$ is determined by a collection of $2n$ numbers $(p_1,q_1,p_3,q_3\dots,p_{2n-1},q_{2n-1})$ (the indices are odd residues modulo $2n$). The definition of these numbers refers to Figure \ref{cross}: 
$$
p_{2i+1}=[A_{2i-3},A_{2i-1},P,R], \ \ q_{2i+1}=[R,Q,A_{2i+3},A_{2i+5}]
$$
where $P=B_{2i-2}\cap B_{2i+2}, R=B_{2i-2}\cap B_{2i+4}, Q=B_{2i}\cap B_{2i+4}$ and 
$[\ ,\ ,\ ,\ ]$ denotes the cross-ratio of four points on a projective line, see  \cite{Sch1,Sch2}. These $2n$ numbers are not independent: they satisfy 8 relations ensuring that the polygon is closed.  Similarly, the dual polygon $L^\ast$ is characterized by the respective cross-ratios $(p^\ast_2,q^\ast_2,\dots,p^\ast_{2n},q^\ast_{2n})$  (the indices are even residues modulo $2n$). It is easy to see that $p^\ast_{2i}=q_{2i-1}$ and $q^\ast_{2i}=p_{2i+1}$. An $n$-gon $L$ is $m$-self-dual if and only if $p_i=p^\ast_{i+m}, q_i=q^\ast_{i+m}$, and hence iff $p_i=q_{i+m-1}, q_i=p_{i+m+1}$. This implies $2m$-periodicity of the sequence of cross-ratios: $p_i=p_{i+2m}, q_i=q_{i+2m}$, cf. Section \ref{m<n}.
}
\end{remark}

\section{The case $m=n$} \label{m=n}

In this section, we consider $n$-gons with odd $n$  and with every vertex dual to the opposite side of a projectively equivalent $n$-gon. According to Proposition \ref{symmetric}, the bilinear form $ F$ is symmetric in this case, so it determines a (complex) Euclidean structure in space $V$. In an 
appropriate coordinate system, the projective duality becomes the {\it polar duality}
$$
(a,b,c)\mapsto \{ax+by+cz=0\}.
$$
(Geometrically, this means that we apply to a point of the standard Euclidean plane 
the inversion in the unit circle centered at 0, then reflect the point in 0, and then take a line through the obtained point perpendicular to the position vector of this point.) For this duality, we will use the notations $A\mapsto A^\perp\mapsto(A^\perp)^\perp=A$.
Two polygons, $n$-self-dual with respect to this duality, are projectively equivalent if and only if they are 
$O(3,\C)$-equivalent.

For a polygon $A_1A_3\dots A_{2n-1}$, construct the star-like polygon $C_1C_2\dots C_n$ where $C_i=A_{1+(i-1)(n-1)}$.

\begin{lemma}\label{star}
The polygon $A_1A_3\dots A_{2n-1}$ is $n$-self-dual $($with respect to the polar duality$)$ if and
only if $C_{i+1}\in C_i^\perp$ for all $i$ $($with $i$  a residue modulo $n)$.
\end{lemma}

\paragraph{\it Proof.} Obvious. \proofend

This leads to a simple explicit construction of all $n$-self-dual $n$-gons. Fix a point $C_1$, then a point 
$C_2\in C_1^\perp$ not equal to $C_1$ (the latter is relevant only if $C_1\in C_1^\perp$). Notice that, 
modulo the action of $O(3,\C)$, there are four choices of the pair $C_1,C_2$ (depending on possible 
incidences $C_1\in C_1^\perp,C_2\in C_2^\perp$). Then choose $C_i\in C_{i-1}^\perp,C_i\ne
C_{i-2},C_{i-1}$ for $i=2,\dots,n-1$, with the additional requirement $C_{n-1}\ne C_1$. In conclusion, we put $C_n=C_{n-1}^\perp\cap C_1^\perp$. After this, we redenote the points, $C_i=A_{1+(i-1)(n-1)}$,
and get an $n$-self-dual $n$-gon $A_1A_3\dots A_{2n-1}$. Moreover, up to a projective equivalence preserving the numeration of vertices, this construction gives all $n$-self-dual $n$-gons, one time each.

In particular, the moduli space of $n$-self-dual $n$-gons has dimension $n-3$ (each of the points 
$C_3,\dots C_{n-1}$ is arbitrarily chosen within a line with finitely many punctures). 

\paragraph{Pentagons.} For an arbitrary $n\ge4$, the moduli space of all $n$-gons has dimension 
$2n-8$ ($2n$ for $n$ vertices, $-8$ for the action of the group $PSL(3,\C)$). In general, this exceeds the  dimension $n-3$ of $n$-self-dual $n$-gons, but for $n=5$ the two numbers coincide: both equal 2. 
Moreover, the following holds.

\begin{proposition}\label{pentagons}
Every pentagon is 5-self-dual.
\end{proposition}

\begin{figure}[hbtp]
\centering
\includegraphics[width=1.8in]{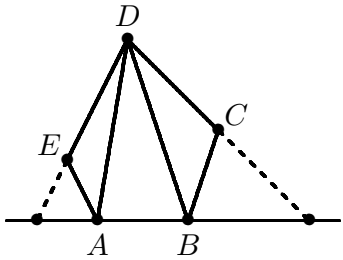}
\caption{Pentagons are self-dual}
\label{pentagon}
\end{figure}

\paragraph{\it Proof.}  We will prove this for a ``generic" pentagon with no three vertices collinear; 
the general case can be resolved by a transition to limit. For a pentagon $ABCDE$ there are 5 cross-ratios $\rho_A,\dots,\rho_E$: $\rho_A$ is defined as the cross-ratio of the lines $AB,AC,AD,AE$, and the other four are defined in a similar way. 
These cross-ratios projectively determine a pentagon: the points $A,B,C,D$  can be moved to  chosen locations, and after that the lines $BE$ and $CE$ are determined 
by $\rho_B$ and $\rho_C$. (Certainly, the five cross-ratios are not independent: generally, two of them 
determine the rest.) For the dual pentagon, these cross-ratos are $\rho_{AB},\dots,\rho_{AE}$ where
$\rho_{AB}$ is the cross ratio of the points $DE\cap AB,A,B,DC\cap AB$ on the line $AB$, and the other four are defined in a similar way. Figure \ref{pentagon} shows that $\rho_E=\rho_{BC}$, and four similar equalities hold as well. \proofend

\paragraph{Poncelet polygons.} We begin with the following easy statement.

\begin{lemma}\label{conic}
Let $C\subset P$ be a non-degenerate conic, and let $L=E_1\dots E_n$ be an $n$-gon inscribed in 
$C$. Then the $n$-gon whose sides are tangent to $C$ at points $E_1,\dots,E_n$
is $($projectively equivalent to the$)$ dual to $L$. More precisely, there exists a projective isomorphism $P\to P^\ast$  that takes $E_i$ to the tangent line to $C$ at $E_i$.
\end{lemma}

\paragraph{\it Proof.} Since all non-degenerate conics are  projectively equivalent, we may assume that $C$ is a unit circle in the Euclidean plane. Let $E'_i$ be the point of $C$ opposite to $E_i$. Then the polygon $L$ is projectively equivalent to the polygon $L'=E'_1\dots E'_n$ and the tangent to $C$ at 
$E'_i$ is polar dual to $E_i$. \proofend 

\begin{figure}[hbtp]
\centering
\includegraphics[width=3in]{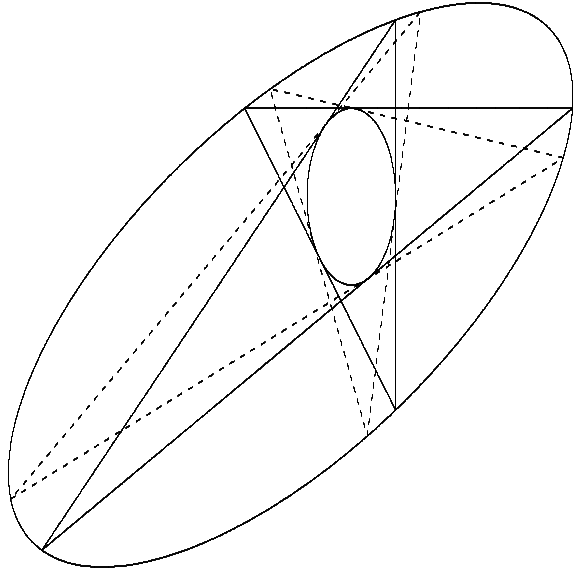}
\caption{Poncelet pentagons}
\label{Poncelet}
\end{figure}

An $n$-gon with an odd $n$ is called a {\it Poncelet polygon} if it is both inscribed into a 
non-degenerate conic and circumscribed about a non-degenerate conic (see Figure \ref{Poncelet}).

\begin{proposition}\label{poncelet}
Every Poncelet $n$-gon is $n$-self-dual.
\end{proposition}

\paragraph{\it Proof.} This follows from Lemma \ref{conic} and the following known result (\cite{Sch,L-T}). Let $A_1A_3\dots A_{2n-1}$ be an $n$-gon (with an odd $n$) inscribed into a conic $C$ and circumscribed about a conic $C'$. Then there exists a projective involution $h\colon P\to P$ such that $h(A_i)$ is the tangency point of $A_{i+n-1}A_{i+n+1}$ and $C'$. \proofend

Two more remarks. The first is that any non-degenerate pentagon is a Poncelet polygon, so 
Proposition \ref{pentagons} follows from Proposition \ref{poncelet}. The second is the following proposition.

\begin{proposition}\label{ponceletmod}
For every odd $n\ge5$, the projective moduli space  of Poncelet $n$-gons is two-dimensional.
\end{proposition}

\paragraph{\it Proof.}  The conics  $C,C'$ from the definition of Poncelet polygons determine a one-parameter family $\mathcal F$ of conics that have four common tangents. Generically, there exists a unique, up to a projective equivalence,  such family $\mathcal F$ (to specify this family, it suffices to fix a generic quadruple of lines). For every $C\in\mathcal F$, there exists a finite number of $C'\in\mathcal F$ such that some $n$-gon inscribed in $C$ is circumsribed about $C'$ (see \cite{G-H} for an explicit condition due to Cayley). Moreover, for such a pair $C,C'$, every point of $C$ is a vertex of such an $n$-gon (Poncelet's theorem, see \cite{B-K-O-R}). Thus, a projective class of a Poncelet $n$-gon is determined by two independent choices: the choice of a $C\in\mathcal F$ and the choice of a point in $C$. \proofend 

Thus, for an odd $n>5$, Poncelet $n$-gons form a small fraction of the space of $n$-self-dual $n$-gons.

\begin{figure}[hbtp]
\centering
\includegraphics[width=5in]{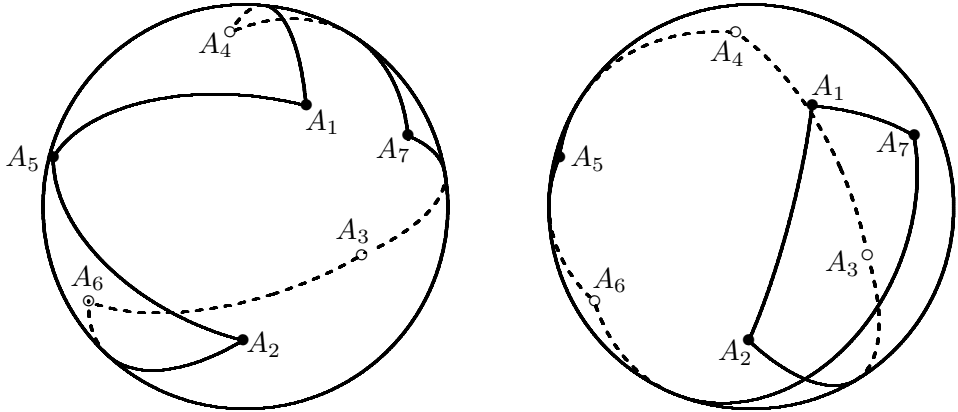}
\caption{Self-dual polygon on the unit sphere}
\label{euclidean}
\end{figure}

\paragraph{The real Euclidean case.} Let $L$ be a real $n$-self-dual $n$-gon. Then $F$ is a real symmetric bilinear form (determined up to real non-zero factors), and there are two possibilities: the form $F$ may be definite or indefinite. In the definite (Euclidean) case, the construction of a self-dual polygon given in the beginning of this section looks especially simple. Consider the unit sphere $S\subset\R^3$. Choose an arbitrary point $C_1\in S$. Then choose a point $C_2$ at the distance  $\displaystyle
\frac\pi2$ from $C_1$. Then choose a point $C_3$ at the distance  $\displaystyle\frac\pi2$ from $C_2$, not equal to $\pm C_1$, then choose $C_4,C_5,\dots$. The last choice will be slightly different from the preceding ones: we choose the point $C_{n-1}$ at the distance  $\displaystyle\frac\pi2$ from $C_{n-2}$, not equal to $\pm C_{n-3}$, and also not equal to $\pm C_1$. After this, we denote by $C_n$ a point at the distance  $\displaystyle\frac\pi2$ from each of the points $C_{n-1}$ and $C_1$. (There are two such points, they form the intersection of two different great circles.) Then we put $C_i=A_{1+(i-1)(n-1)}$ and project the polygon $A_1A_3\dots A_{2n-1}$ onto $P$. This is our self-dual polygon (see Figure \ref{euclidean}).

 It is natural to ask, which $n$-self-dual $n$-gons correspond to definite forms. A partial answer to this
 question is provided by the following proposition. 
 
 \begin{proposition}\label{convex}
 If a real $n$-self-dual $n$-gon is projectively equivalent to an affine convex $n$-gon, then the corresponding symmetric bilinear form is definite.
 \end{proposition}
 
\begin{figure}[hbtp]
\centering
\includegraphics[width=2in]{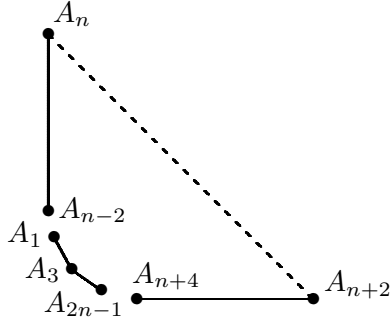}
\caption{Convex septagon: proof of Proposition \ref{convex}}
\label{conv7}
\end{figure}

 \paragraph{\it Proof.} Let $A_1A_3\dots A_{2n-1}$ be our convex polygon. We will use only the convexity of the heptagon $A_1A_3A_{n-2}A_nA_{n+2}A_{n+4}A_{2n-1}$.  We will assume that the points $A_{n-2},A_n,A_{n+2},A_{n+4}$ have projective coordinates $(0:1:1),(0:1:0),(1:0:0),(1:0:1)$.
Then the lines $B_{n-1}=A_{n-2}A_n,\ B_{n+1}=A_nA_{n+2},\ B_{n+3}=A_{n+2}A_{n+4}$ are, respectively, the $y$-axis, the line at infinity, and the $x$-axis (see Figure \ref{conv7}). Therefore $B^\ast_{n-1}=(1,0,0),B^\ast_{n+1}=(0,0,1), B^\ast_{n+3}=(0,1,0)$. The matrix of the isomorphism $G=F^{-1}\colon P^\ast\to P$ which takes $B^\ast_j$ into $A_{j+n}$ is symmetric by Proposition \ref{symmetric}. 
Let it be \[G=\left[\begin{array} {ccc} a&b&c\\ b&d&e\\ c&e&f\end{array}\right].\] Then $A_{2n-1}=
G(B^\ast_{n-1})=(a,b,c), A_1=G(B^\ast_{n+1})=(c,e,f), A_3=G(B^\ast_{n+3})=(b,d,e)$, and the affine coordinates of the points $A_{2n-1},A_1,A_3$ are$$x_{2n-1}=\frac ac , y_{2n-1}=\frac bc,\ x_1=\frac cf,\ 
y_1=\frac ef,\ x_3=\frac be,\ y_3=\frac de.$$ 

The conditions of convexity of our heptagon are $0<x_3<x_1<x_{2n-1}<1,$ $1>y_3>y_1>y_{2n-1}>0$ and \[\det\left[\begin{array} {cc} 
x_1-x_{2n-1}&y_1-y_{2n-1}\\ x_1-x_3&y_1-y_3\end{array}\right]>0;\]the latter means$$x_1y_{2n-1}+x_3y_1+x_{2n-1}y_3-x_1y_3-x_3y_{2n-1}-x_{2n-1}y_1>0,$$ or$$2\cdot\frac bf+\frac ac\cdot\frac de-\left(\frac cf\cdot\frac de+\frac be\cdot bc+\frac ac\cdot\frac ef\right)>0.$$After multiplication by the positive number $c^2f^2x_{2n-1}y_1=acef$, this inequality becomes$$2acbe+a^2df-(ac^2d+ab^2f+
a^2e^2)>0,$$that is, $a\cdot\det G>0$. Also it follows from the convexity inequalities that $ad-b^2=f^2x_1y_1(x_{2n-1}y_3-x_3y_{2n-1})>0$. These two inequalities show that the form $F$ is 
definite. \proofend

It should be noted that, as it is seen from the construction above (with a sequence of points on the sphere), a self-dual polygon with a definite form does not need to be convex. Still, it is true that a 
pentagon has a definite form if and only if it is projectively equivalent to a convex pentagon (we leave
a proof to the reader).

\section{The case $m<n$} \label{m<n}

Let $L=A_1A_3\dots A_{2n-1}$ be an $m$-self-dual polygon with $m<n$. By Proposition
 \ref{symmetric}, the form $F$ in this case is not symmetric. If $B$ is a point or a line in $P$, then
$B^\perp$ is defined as $\{y\in P\mid F(y,x)=0\ \mbox{for all}\ x\in B\}$. There arises a projective transformation $G\colon P\to P,\ G(B)=(B^\perp)^\perp$, and obviously, in term of matrices, 
$G=F^{-1}F^t$. 

\begin{lemma} \label{shift}
If $L$ is $m$-self-dual then $G(A_i)=A_{i+2m}$.
\end{lemma}

\paragraph{\it Proof.} 
One has $A_i^\perp=A_{i+m-1}A_{i+m+1}$, and hence
$G(A_i)=(A_i^\perp)^\perp=A_{i+m-1}^\perp\cap A_{i+m+1}^\perp=(A_{i+2m-2}A_{i+2m})\cap(A_{i+2m}A_{i+2m+2})=A_{i+2m}.$
 \proofend

Thus, $G$ makes a non-trivial cyclic permutation of vertices of $L$, and, in particular, $G^r=\id$ where 
$r=\displaystyle\frac n{(m,n)}$. The following result is an elementary fact from linear algebra.

\begin{lemma}\label{canonical}
Let $F$ be a non-degenerate non-symmetic bilinear form in $V$. Then there exists a basis in $V$ with respect to which $F$ has one of the following matrices:$$H_\varphi=\left[\begin{array} {ccc} \phantom{-}\cos\varphi&\sin\varphi&0\\ -\sin\varphi&\cos\varphi&0\\ \phantom{-} 0&0&1\end{array}\right],\  J=\left[\begin{array} {ccc} \phantom{-} 1&1&0\\ -1&0&0\\ \phantom{-} 0&0&1\end{array}\right],\ K=\left[\begin{array} {ccc} \phantom{-} 1&1&0\\ -1&0&1\\ \phantom{-} 0&1&0\end{array}\right].$$
\end{lemma}

\paragraph{\it Proof.} There is a unique decomposition $F=F_++F_-$ where $F_+$ is symmetric and $F_-$ is skew-symmetric. Let $W=\Ker F_-$; since $F_-\ne0$, $\dim W=1$.\smallskip

\noindent{\it Case} 1: $F_+|_W\ne 0,\rank F_+=3$. Let $e_3\in W,\ F_+(e_3,e_3)=1$. Let $Z$ be the orthogonal complement to $W$ with respect to $F_+$. Choose $e'_1,e'_2\in Z$ with $F_+(e'_i,e'_j)=\delta_{ij}$. Then choose a rescaling $e_1=\alpha e'_1, e_2=\alpha e'_2$ such that $F_+(e_1,e_1)^2+
F_-(e_1,e_2)^2=1$. (Notice that $F_+(e'_1,e'_1)^2+F_-(e'_1,e'_2)^2=\det F\ne0$.) Then the matrix of 
$F$ with respect to the basis $e_1,e_2,e_3$ is $H_\varphi$ with some (complex) 
 $\varphi\ne\displaystyle\frac{k\pi}2$.\smallskip

\noindent{\it Case} 2: $F_+|_W\ne 0, \rank F_+=2$.  Let $e_3, Z$ denote the same as in Case 1, let 
$0\ne e'_2\in\Ker F_+$, and let $e_1\in Z-\Ker F_+,\ F_+(e_1,e_1)=1$. Choose $\alpha$ such that, for 
$e_2=\alpha e'_2,\ F_-(e_1,e_2)=1$. Then the matrix of $F$ with respect to the basis $e_1,e_2,e_3$ is
$J$.\smallskip

\noindent{\it Case} 3: $F_+|_W\ne 0, \rank F_+=1$. Let $e_3$ denote the same as in Cases 1 and 2,
and let $e_1,e_2$ be a basis in $\Ker F_+$ such that $F_-(e_1,e_2)=1$. Then the matrix of $F$ with respect to the basis $e_1,e_2,e_3$ is $\displaystyle{\left[\begin{array} {ccc} \phantom{-} 0&1&0\\ 
-1&0&0\\ \phantom{-} 0&0&1\end{array}\right]}$ which is $H_{\frac\pi2}$.\smallskip

\noindent{\it Case} 4: $F_+|_W=0, \rank F_+=3$. Choose a non-zero vector $e'_3\in W$. Let $Z$ be the orthogonal complement to $W$, and let $C$ be the ``light cone" $\{x\in V\mid F_+(x,x)=0\}$. Since $C\not\subset Z$, we can choose an $e'_2\in V$ such that $F_+(e'_2,e'_2)=0,\ F_+(e'_2,e'_3)=1$. Let 
$U\in V$ be the subspace spanned by $e'_2$ and $Y$ be the orthogonal complement of $U$. The 
intersection $Y\cap Z$ is not contained in $W+U$: if a linear combination of $e'_2$ and $e'_3$ is orthogonal to both $e'_2$ and $e'_3$, then it must be 0. Take  $e_1\in Y\cap Z$ with 
$F_+(e_1,e_1)=1$. Then $F_-(e_1,e'_2)\ne 0$ (otherwise $F_-$ would have been zero). Put
$e_2=\alpha e'_2,\ e_3=\alpha^{-1} e'_3$ in such a way that $F_-(e_1,e_2)=1$.  Then the matrix of $F$ with respect to the basis $e_1,e_2,e_3$ is $K$.\smallskip  

\noindent{\it Case} 5: $F_+|_W=0, \rank F_+<3$. Take a 1-dimensional space $U\subset\Ker F_+$. If 
$U=W$, then $\Ker F\supset W$ is non-zero, so $F$ is degenerate. If $U\ne W$, then both $F_+,F_-$ are zero on $U\oplus W$, which also means that $F$ is degenerate. \proofend

Return now to the bilinear form $F$ related to our $m$-self-dual $n$-gon $L$.

\begin{proposition}\label{rotation}   
In an appropriate coordinate system, the matrix of $F$ is $H_\varphi$ with $r\varphi\in\pi\mathbb Z$. Moreover, if the $n$-gon $L$ is simple, then $s\varphi\notin\pi\mathbb Z$ for any positive $s<r$.
\end{proposition}

\paragraph{\it Proof.} According to Lemma \ref{canonical}, the matrix of $F$, in an appropriate
basis, is $H_\varphi,\ J$, or $K$. But
$$J^{-1}J^t=\left[\begin{array} {ccc} -1& \phantom{-}0&0\\  
\phantom{-}2&-1&0\\  \phantom{-}0& \phantom{-}0&1\end{array}\right],\ K^{-1}K^t=\left[\begin{array} {ccc}  1&-2&0\\ 0&\phantom{-}1&0\\ 2&-2&1\end{array}\right]$$
and neither of these two matrices has  finite order (for both, the Jordan form contains a non-trivial Jordan block). On the other hand, the matrix of 
$G$ is $(H_\varphi)^{-1}H_\varphi^t=H_{-\varphi}^2=H_{-2\varphi}$, that is (again, in an appropriate coordinate system), $H_{-2\varphi}A_i=A_{i+2m}$. First, this shows that $H_{-2\varphi}^r=I$, that is, 
$2r\varphi$ is a multiple of $2\pi$. Second, if $2s\varphi$ is a multiple of $2\pi$ for a positive $s<r$, 
then $A_i=A_{i+2sm}$ where $sm$ is not a multiple of $n$, so our polygon is not simple. \proofend

Thus, $L$ contains $(m,n)$ regular $r$-gons, $A_iA_{i+2m}\dots 
A_{i+(r-1)m},\ i=1,2,$ $\dots,(m,n)$. [By a regular $n$-gon we understand a (maybe, self-intersecting) 
$n$-gon in the Euclidean plane with all lengths of the sides equal and all angles equal; thus, there are two projectively different types of regular pentagons, three projectively different types of regular heptagons, and so on.]  If $(m,n)=1$, then $L$ itself is regular (so, in this case, an $m$-self-dual $n$-gon is projectively unique). If $(m,n)>1$, then this uniqueness, in general, does not hold. Below, we give an explicit construction of all $m$-self-dual $n$-gons which will demonstrate this non-uniqueness. 

First, notice that, in our case, the projective duality has a simple geometric description: we consider a Euclidean plane (with a fixed origin) and, for a point $A\ (\ne0)$, the dual line $A^\ast$ is obtained from 
the polar dual $A^\perp$ by a clockwise rotation about the origin by the angle of $\displaystyle
\frac{\pi(n-m)}n$.

Now, let us construct an arbitrary $m$-self-dual $n$-gon. In addition to $r=\displaystyle\frac n{(m,n)}$, put $k=\displaystyle\frac m{(m,n)}$ and also $d=(m,n)$; thus, $mr=kn$. In the Euclidean plane with a fixed origin $O$, choose an arbitrary point $A_1$. Then successive counter-clockwise rotations by the angle $\displaystyle\frac{2\pi m}n$ give the points $A_{2m+1}, A_{4m+1},\dots, A_{2(r-1)m+1}$, and also the lines $A^\ast_1=A_mA_{m+2},A^\ast_{2m+1}=A_{3m}A_{3m+2},A^\ast_{4m+1}=A_{5m}A_{5m+2},\dots, A^\ast_{2(r-1)m+1)}=A_{(2r-1)m}A_{(2r-1)m+2}$.  Of the numbers $m,3m,5m,\dots,2(r-1)m$ one
is $d$ modulo $2n$ ($d=um+v\cdot2n$ for a unique $u,\ 0\le u<2r$, and this $u$ must be odd). So, one of our lines should be $A_dA_{d+2}$; choose a point $A_d$ on this line. This choice gives also the points $A_{2m+d},A_{4m+d},\dots,A_{2(r-1)m+d}$ and the lines $A_{d+m-1}A_{d+m+1},
A_{d+3m-1}A_{d+3m+1},$ $A_{d+5m-1}A_{d+5m+1},\dots,A_{d+(2r-1)m-1}A_{d+(2r-1)m+1}$. By the way, one of these lines will be $A_{2n-1}A_1$. Our next choice will be a point $A_{d+2}$, again on the line $A_dA_{d+2}$. This will give us $r$ additional points (including $A_{d+2}$) and $r$ lines, dual to these points. One of these lines will be $A_1A_3$, and we choose a point $A_3$ on it. One of the lines coming with this point will be $A_{d+2}A_{d+4}$, and we choose a point $A_{d+4}$, and so on. Proceeding in these way, we choose the points in the following order: $A_1,A_d, A_{d+2},A_3,A_{d+4},
A_5,A_{d+6}A_7,\dots ,A_{2d-3},A_{d-2}$. Here we stop: the next choice should be $A_{2d-1}$, but this
point will appear as the intersection of the line $A_{2d-1}A_{2d+1}$ coming with the point $A_d$ and the line $A_{2d-3}A_{2d-1}$ coming with the point $A_{d-2}$. After that, we have the points $A_1,A_3,
A_5,\dots,A_{2d-1}$, and hence we have all the vertices of our polygon.

The projective symmetry $G,\ G(A_i)=A_{i+2m}$, shows that every $m$-self-dual $n$-gon is also 
$m'$-self-dual for every $m'\equiv em\bmod n$ where $e$ is odd. In particular, if $(e,r)=1$, then 
$m$-self-dual $n$-gons and $m'$-self-dual $n$-gons are the same $n$-gons (although their 
self-dualities involve different projective isomorphisms $P\to P^\ast$). On the other hand, if $n$ is odd, then we see that every $m$-self-dual $n$-gon is also $n$-self-dual, that is, it belongs to the class of polygons considered in Section 3.

\begin{figure}[hbtp]
\centering
\includegraphics[width=1.4in]{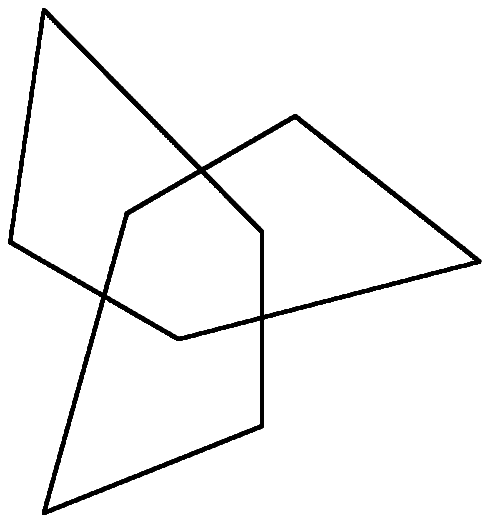}
\qquad
\includegraphics[width=1.4in]{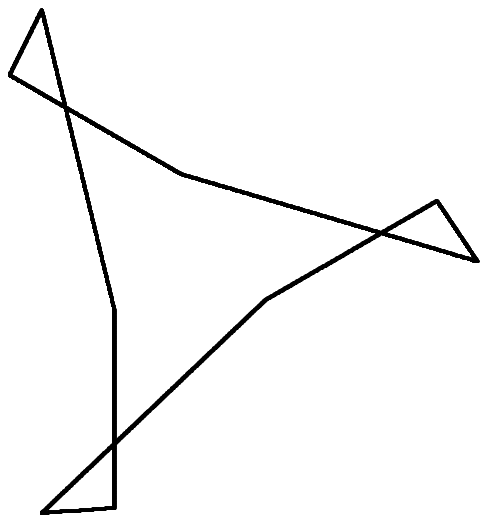}
\caption{Two 3-self-dual nonagons}
\label{9-gons}
\end{figure}

\begin{figure}[hbtp]
\centering
\includegraphics[width=2.2in]{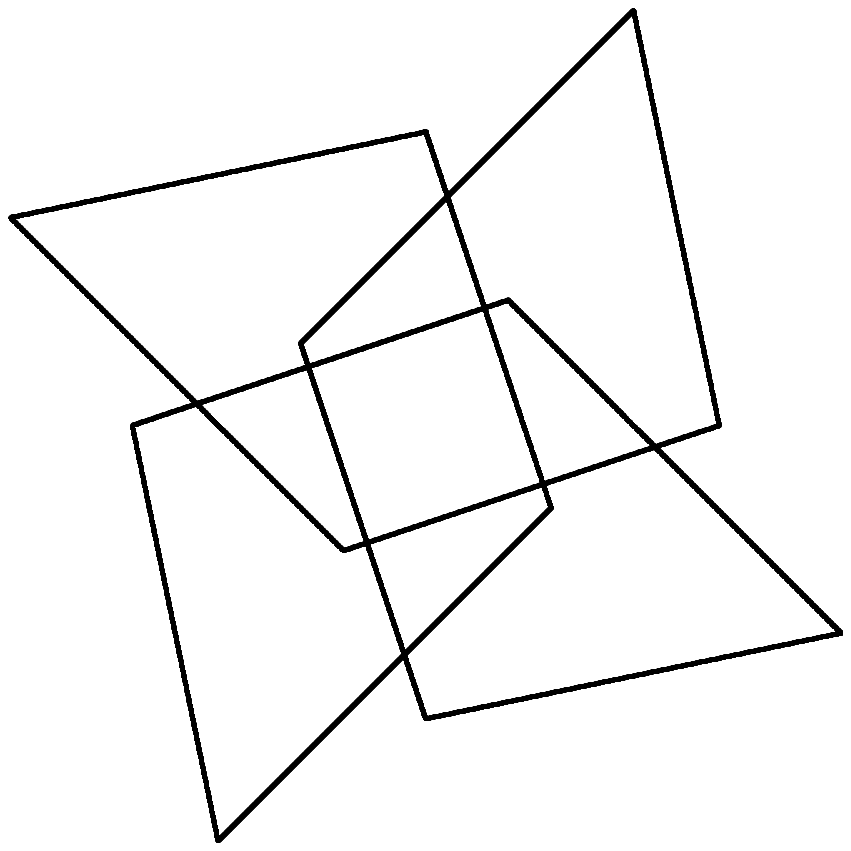}
\qquad
\includegraphics[width=2.2in]{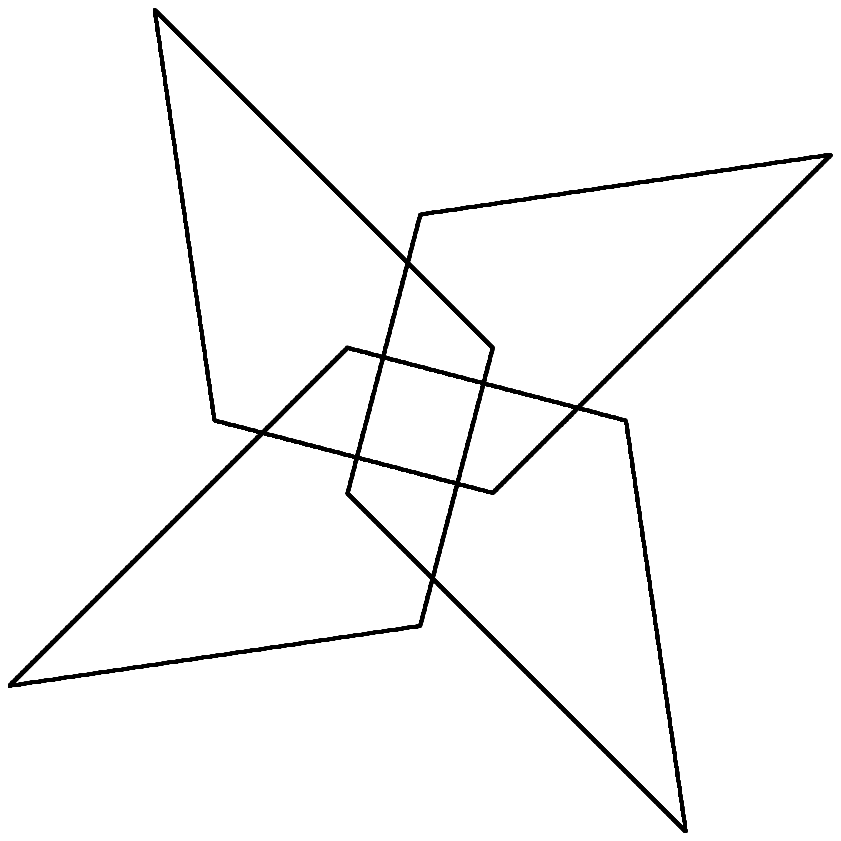}
\caption{Two 3-self-dual dodecagons}
\label{12-gons}
\end{figure}

Figure \ref{9-gons} shows two 3-self-dual 9-gons. As was remarked above, they are also 9-self-dual (with respect to the polar duality). Note that no $n$-gon with $n$ even is $n$-self-dual (the definition of $m$-self-duality requires that $m$ is odd); but they must be centrally symmetric (with respect to the affine chart considered in this section), and the 12-gons of Figure   \ref{12-gons} are   centrally symmetric indeed.                                                                                                                                                                                                                                                                                                                                                                                  

\begin{proposition}\label{dimension}
Let $m<n$ and $(m,n)>1$. Then the moduli space of $m$-self-dual $n$-gons has the dimension 
$(m,n)-1$, for $n\ne2m$, and $(m,n)-3$, for $n=2m$.
\end{proposition}

\paragraph{\it Proof.} If the basis, in which the bilinear form $F$ has the canonical form $H_\varphi$, is chosen then, to specify our $m$-self dual $n$-gon (with $m<n$), we need to choose a point $A_1$ in the plane (minus one point), which depends on two parameters, and then $d-2,\ d=(m,n)$, points $A_d,A_{d+2},
A_3,A_{d+4},A_5,\dots,A_{d-2}$ which provide $(m,n)-2$ more parameters, with the total of $(m,n)$. From this number, we need to subtract the number of parameters on which the basis for a given form $F$ depends (in other words, the dimension of the Lie group of linear transformations of $V$ that preserve the form $H_\varphi$). If $2\varphi$ is not a multiple of $\pi$, which corresponds to the case $n\ne2m$, this dimension is 1. This is seen from Case 1 of the proof of Lemma \ref{canonical}: the choice of $e_3$ provides no parameters (it is two-valued), then we choose $e'_1$ and $e'_2$ on the conic 
$F_+(x,x)=0$ with the condition $F_+(e'_1,e'_2)=0$, which provides one parameter, and then we multiply both $e'_1$ and $e'_2$ by the same complex number which we determine from a quadratic equation. So the total number of free parameters in this case is 1, and the dimension of the moduli space is $(m,n)-1$. If $n=2m$, then $\cos\varphi=0$, and this is Case 3 of the proof of Lemma 
\ref{canonical}. In this case, the choice of $e_3$ does not provide any parameters while $e_1$ and $e_2$ are chosen up to the action of the group of transformations of $Z$ preserving the form $F_-$; it is $SL(2;\C)$, the dimension is 3. Thus, if $n=2m$, then the dimension of the moduli space is 
$(m,n)-3=m-3$. \proofend

Notice in conclusion that our results show that the space of moduli of $m$-self-dual hexagons has dimension 0, whatever $m$ is. Actually, the only self-dual hexagon (for any $m$) is the regular hexagon; again, we leave the details to the reader.

\section{Polygonal curves} \label{polycurves}

A {\it polygonal curve} is a real polygon $A_1A_3\dots A_{2n-1}$ with the following two (independent) additional structures. (1) For every even $i$, one of the two segments into which the points 
$A_{i-1},A_{i+1}$ cut the real projective line $B_i$ is chosen; we will refer to this segment as  an {\it edge} of the polygonal curve. (2) For every odd $i$, one of the two pairs of vertical angles formed by the real projective lines $B_{i-1},B_{i+1}$ is fixed; we will refer to these angles as {\it exterior angles} of
the polygonal curve. Thus, every real $n$-gon gives rise to $2^{2n}$ polygonal curves.

For a polygonal curve $A_1A_3\dots A_{2n-1}$ (with its additional structures), there arises a {\it dual polygonal curve} in dual real projective plane $P^\ast$. This is the polygon $B^\ast_2B^\ast_4\dots B^\ast_{2n}$ with the following edges and exterior angles. The edge $B^\ast_{2i}B^\ast_{2i+2}$ is formed by points of $P^\ast$ dual to the lines through $A_{2i+1}$ in $P$ contained in the exterior angles at $A_{2i+1}$ of the given polygonal curve. The exterior angles at $B^\ast_{2i}$ is formed by lines in 
$P^\ast$  dual to the points of the edge $A_{2i-1}A_{2i+1}$ of the given polygonal curve.

Let $m\le n$ be a positive odd number. A polygonal curve $A_1A_3\dots A_{2n-1}$ is called 
$m$-self-dual if there exists a projective isomorphism $P\to P^\ast$ which takes $A_i$ into 
$B^\ast_{i+m}$ and also takes edges and exterior angles of the polygonal curve $A_1A_3\dots 
A_{2n-1}$ into edges and exterior angles of the polygonal curve $B^\ast_2B^\ast_4\dots B^\ast_{2n}$.

\begin{proposition} \label{choices}
Let $L=A_1A_3\dots A_{2n-1}$ be an $m$-self-dual $n$-gon. Then, of the $2^{2n}$ polygonal curves arising from the polygon $L$, $2^{(m,n)}$ are $m$-self-dual.
\end{proposition}

\paragraph{\it Proof.} Choose edges $A_1A_3,A_3A_5,\dots A_{2d-1}A_{2d+1}$ in an arbitrary way. If $m<n$, then apply to these edges $r=\displaystyle\frac nd$ consecutive rotations by the angle 
$\displaystyle\frac {2\pi m}n$; we will get a full set of edges. Then apply the duality to these edges, and
this gives a choice of exterior angles for the dual polygon $L^\ast$. The projective isomorphism between $L$ and $L^\ast$ makes these angles exterior angles for $L$, and $L$ becomes an $m$-self-dual 
polygonal curve. Obviously, this construction gives all $m$-self-dual $n$-gonal curves. \proofend

Polygonal curves are  natural polygonal counterparts to smooth curves. Let us describe a unifying point of view. 

\paragraph{The space of contact elements and projective duality.}
A contact element of the real projective plane is a pair $(A,B)$ where $A\in P$ is a point, $B\subset P$ is a line, and  $A\in B$. Denote the space of contact elements by $F$ (it is naturally identified with the space of full flags in $\R^3$). One has two projections $\pi_1\colon F\to P$ and $\pi_2\colon F\to P^\ast$ defined by the formulas: $\pi_1(A,B)=A,\ \pi_2(A,B)=B$. The space $F$ has a contact structure (a non-integrable two-dimensional distribution) defined by the condition that the velocity of point $A$ lies in the line $B$. The fibers of the projections $\pi_{1,2}$ are Legendrian curves (curves tangent to the contact distribution). The space of contact elements of the dual plane $P^\ast$ is canonically identified with $F$.
 
Projective duality is easily described in terms of the space of contact elements. Let $\gamma\subset P$ be a smooth curve. Assigning the tangent line to each point of $\gamma$ gives a lift $\Gamma\subset F$; this lifted curve is Legendrian. The curve $\pi_2(\Gamma)$ is the dual curve $\gamma^\ast\subset P^\ast$, and the lift of $\gamma^\ast$ to $F$ is again $\Gamma$. A wave front in $P$ is defined as the $\pi_1$-projection of a smooth Legendrian curve $\Gamma$ in $F$ to $P$; it has singularities (generically, semi-cubical cusps) at the points where $\Gamma$ is tangent to the fibers of $\pi_1$. The dual wave front is $\pi_2(\Gamma)$.

Likewise for polygonal curves. The set of edges of a polygonal curve provides a closed curve in $P$, but its lift to $F$ consists of disjoint arcs of the fibers of $\pi_2$; to connect these arcs by segments of the fibers of $\pi_1$, we need to choose exterior angles. Thus, for an $n$-gonal curve $L$ in $P$, its lift to $F$ is a 
$2n$-gonal curve, whose sides are segments of the  fibers of the alternating projections $\pi_1$ and $\pi_2$, and  whose projection to $P^\ast$ is the dual polygonal curve $L^\ast$.

\section{Self-dual curves} \label{sdcurves}

Many a result from the preceding sections extends to self-dual wave fronts in $\RP^2$. Here we do not 
attempt to give a complete classification of such fronts; instead we describe several classes of examples, including explicit formulas for self-dual curves and fronts. Self-dual curves will be described as Legendrian curves in certain three-dimensional contact manifolds satisfying certain monodromy conditions. 

In this section, $P=\RP^2$ and $V=\R^3$. 
Let $\gamma(t)\subset P$ be a self-dual parameterized closed curve (possibly, with cusps); we assume that the parameter $t$ takes values in $S^1=\R/2\pi\Z$.  The projectively dual curve  also has a parameterization, 
$\gamma^\ast(t)$: the covector $\gamma^\ast(t)$, defined up to a non-zero multiplier, vanishes on the vectors $\gamma(t)$ and $\gamma'(t)$.  As in Section \ref{poldual}, we have a linear isomorphism $f\colon V \to V^\ast$ that takes the line $\gamma(t)\subset V$ to the line $\gamma^\ast(\varphi^{-1}(t))\subset V^\ast$ where $\varphi$ is a diffeomorphism of $S^1$. This diffeomorphism plays the role of the cyclic shift through $m$ in the definition of $m$-self-dual polygons. As in Section \ref{poldual}, we consider the corresponding projective isomorphism $\hat f\colon P\to P^\ast$ and  the bilinear form $F$ on $V$, $F(v,w)=\langle f(v),w\rangle$. As in Section \ref{m<n}, we consider the projective transformation $G=(\hat f)^{-1}\hat f^\ast \colon P\to P,\ G(B)=(B^\perp)^\perp$.

The next lemma is an analog of Lemma \ref{shift}.

\begin{lemma} \label{curveshift}
One has: $G(\gamma(t))=\gamma(\varphi^2(t))$ for all $t$.
\end{lemma}

\paragraph{\it Proof.} The proof is essentially the same as that of Lemma \ref{shift}. First note that the line 
$\{y\in P\ |\ \langle y,\gamma(t)\rangle=0\}$ is the tangent line $T_{\gamma^\ast(t)} \gamma^\ast$ to $\gamma^\ast$ at point $\gamma^\ast(t)$. It follows that $\gamma(t)^\perp = (\hat f)^{-1} (T_{\gamma^\ast(t)} \gamma^\ast)=T_{\gamma(\varphi(t))} \gamma$. Likewise, the point $\{y\in P\ |\ \langle y,T_{\gamma(\varphi(t))} \gamma\rangle=0\}$ is $\gamma^\ast(\varphi(t))$, and therefore $(\gamma(t)^\perp)^\perp = \gamma(\varphi^{2}(t))$.
\proofend

Analogs of the polygons considered in Section \ref{m=n} are the curves for which the diffeomorphism $\varphi$  is an involution. Just like a polygon, a curve may be a multiple of another curve. The next proposition is an
analog of Proposition \ref{symmetric}. 

\begin{proposition} \label{curvesymm}
Assume that a self-dual curve $\gamma$ is not a multiple of another curve. Then $\varphi^2=\id$ if and only if the bilinear form $F$ is symmetric.
\end{proposition}

\paragraph{\it Proof.} $F$ is symmetric if and only if $G$ is the identity. If $\varphi^2=\id$ then, by Lemma \ref{curveshift}, $G(\gamma(t))=\gamma(t)$ for all $t$. Since $\gamma$ contains four points in general position, $G=\id$. Conversely, if $G=\id$ then $\gamma(t)=\gamma(\varphi^2(t))$ for all $t$. Since $\gamma$ is not  a multiple of another curve, $\varphi^2=\id$.
\proofend

An analog of Proposition \ref{convex} holds as well.

\begin{proposition} \label{curvepos}
If $\gamma$ is  a convex self-dual curve such that $\varphi^2=\id$ then the symmetric bilinear form $F$ is definite.
\end{proposition}

\paragraph{\it Proof.} Assume not. Then, in an appropriate coordinate system, the respective quadratic form is $x^2+y^2-z^2$. The light cone $x^2+y^2=z^2$ projects to a circle $C\subset \RP^2$. If  $p\in C$ then the line $p^\perp$ is the tangent line to $C$ at point $p$, and if a point $p$ is inside $C$ then the line $p^\perp$ lies in the exterior of $C$. 

Due to the convexity, $\gamma(t)\notin \gamma(t)^\perp$ for all $t$, hence $\gamma$ does not intersect $C$. Therefore $\gamma$ lies either inside $C$ or outside of it. In the former case, the envelope $\gamma^\ast$ of the lines $\gamma(t)^\perp$ lies outside of $C$ and cannot coincide with $\gamma$. In the later case, one can find a tangent line $\ell$ to $\gamma$, disjoint from $C$; then the point $\ell^\perp\in\gamma^\ast$ lies inside $C$, and again $\gamma$ fails to coincide with $\gamma^\ast$.
\proofend

\paragraph{Spherical curves.}
 Let us consider the case when the symmetric bilinear form $F$ is (positive) definite. Then the correspondence between points and the dual lines is that between pairs of antipodal poles and the corresponding equators on the unit sphere (that doubly covers $\RP^2$). Thus the projective duality moves every point of a curve  $\gamma\subset S^2$ distance $\displaystyle\frac\pi2$ in the normal direction to $\gamma$.

An example of a self-dual curve on $S^2$ is a circle of radius $\displaystyle\frac\pi4$. This circle is included into the family of curves of constant width $\displaystyle\frac\pi2$ that are all self-dual (of course, all distances are in the spherical metric). We interpret curves of constant width as Legendrian curves. 

Let $M$  be the space of oriented geodesic segments of length  $\displaystyle\frac\pi2$ on $S^2$. Then $M$ is diffeomorphic to $SO(3)\cong \RP^3$. Define a two-dimensional distribution $E$ on $M$ by the condition that the velocities of the end points of a segment are perpendicular to the segment. Assign an oriented contact element to a geodesic segment $AB$: the foot point is the midpoint of $AB$ and the direction is the oriented normal to $AB$. This provides an identification of $M$ with the space of oriented contact elements $F$ of $S^2$.

\begin{lemma}
Under the diffeomorphism $F\cong M$, the standard contact structure in $F$ is identified with the distribution $E$.
\end{lemma} 

\paragraph{\it Proof.} The space $E$ is generated by  two vector fields corresponding to the following motions of a geodesic segment $AB$: the rotation of $AB$ about its midpoint, and the rotation of $AB$ about the axis $AB$ (so that the end points are fixed). The velocities of the corresponding motions of the midpoint of $AB$ are orthogonal to $AB$ which proves the lemma. 
\proofend

Thus a curve of constant width can be constructed as a smooth Legendrian curve $A(t) B(t) \subset M,\ t\in \R$ satisfying the monodromy condition $A(\pi)=B(0), B(\pi)=A(0)$. Clearly, there is an abundance of such curves, in particular, analytic ones. This construction gives curves with cusps and inflections as well.   
(A similar approach is used in \cite{B-Zh} to construct billiard tables that possess one-parameter families of periodic trajectories, the case of two-periodic trajectories being that of curves of constant width.)

To construct such a curve, take a closed wave front on the sphere with an odd number of cusps (say, an odd-cusped hypocycloid), place a geodesic segment of length $\displaystyle\frac\pi2$ orthogonally to the front, so that its midpoint is on the front, and use the front as a guide to move the geodesic segment all the way around until its end points swap their positions.

\begin{figure}[hbtp]
\centering
\includegraphics[width=2in]{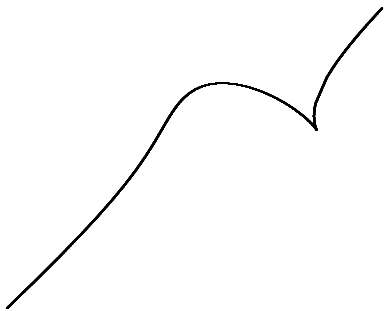}
\caption{Self-dual cubic curve with one inflection and one cusp}
\label{cubic}
\end{figure}

Curves of constant width $\displaystyle\frac\pi2$ in $S^2$ project to $\RP^2$ as contractible curves. 
Similarly one can construct self-dual non-contractible curves in $\RP^2$. For this, one needs to modify the above monodromy condition: $A(\pi)=B(0), B(\pi)=-A(0)$. A  non-contractible curve necessarily has an odd number of inflections, and therefore, by self-similarity, an odd number of cusps. An example is a cubic curve $y=x^3$; after a projective transformation, this curve looks like shown on Figure \ref{cubic}.

\paragraph{Rotationally symmetric curves.} Let us consider analogs of the polygons studied in Section \ref{m<n}. Assume that $\gamma$ is a self-dual curve that lies in the affine plane (i.e., is disjoint from the line at infinity) and is star-shaped with respect to the origin $O$ (i.e., the tangent lines to $\gamma$ do not pass through $O$). We assume that $\gamma$ is not a multiple of another curve.  We allow $\gamma$ to have inflections and cusps.  Assume that the bilinear form $F$ is not symmetric, so, by Proposition \ref{curvesymm}, $\varphi^2\neq \id$.

Arguing as in Section \ref{m<n}, we choose a coordinate system in which the mapping $G$ is a rotation through some angle $\alpha$ (the cases of Jordan blocks and complex angles in Lemma \ref{canonical} are excluded because in these cases orbits of $G$ would have accumulation points at infinity, in contradiction with Lemma \ref{curveshift}). If $\alpha$ is $\pi$-irrational then the orbits of $G$ are dense in a circle and, by Lemma \ref{curveshift}, $\gamma$ is a circle. Otherwise, $\alpha=\displaystyle\frac{2\pi p}{q}$ where $p$ and $q$ are co-prime. Thus $q$ is the least period of $G$, and hence of the circle diffeomorphism $\varphi^2$ as well. Choose a parameterization $\gamma(t)$ so that $\displaystyle\varphi^2(t)=t+\frac{2\pi r}{q}$ where $r$ and $q$ are also co-prime. To summarize, we have the following analog of Lemma \ref{shift}.

\begin{lemma} \label{cshift}
One has: $\displaystyle G(\gamma(t))=\gamma\left(t+\frac{2\pi r}{q}\right)$.
\end{lemma}

Note that the rotation number of $\gamma$ about the origin $O$ equals the least positive $k$ such that $kr=p$ mod $q$. 

Let us describe an explicit construction of such self-dual curves. Let $H$ be the rotation about the origin through angle $\displaystyle\frac{\pi p}{q}$ (so that $H^2=G$) and set $c=\displaystyle\frac{\pi r}{q}$. 
As in Section \ref{m<n}, the projective duality has a simple geometric description in terms of the Euclidean metric: for point $A$, the dual line $A^\ast$ passes through the point $H(A)/|A|^2$ and is orthogonal to the vector $H(A)$. 

Self-dual curves again can be described as Legendrian curves. Let $M$ consists of pairs of vectors $(u,v)$ such that $H(u)\cdot v=1$. Define a contact structure on $M$ by the condition $H(u)\cdot v'=0$ (or, equivalently, $H(u')\cdot v=0$); we leave it to the reader to check that this is indeed a contact structure. Let $(u(t),v(t))$ be a Legendrian curve in $M$ satisfying the monodromy condition $u(t+c)=v(t), v(t+c)=G(u(t))$. Then we can set: $\gamma(t)=u(t)$. The condition $H(u)\cdot v=1$ implies that $v$ belongs to the line $u^\ast$, and the Legendrian condition $H(u)\cdot v'=0$ that this line is tangent to the curve $\gamma$ at point $\gamma(t+c)$. Thus $\gamma$ is  self-dual.

Now we give explicit formulas. 

\begin{proposition} \label{expform}
Let $\beta(t)$ be a smooth function such that $|\beta(t)|<\displaystyle\frac\pi4$ and $\beta(t+c)=-\beta(t)$. Let $\rho_1(t)$  and $\rho_2(t)$ satisfy the differential equations 
\begin{equation} \label{form}
\rho'_1=(\beta'+1) \tan 2\beta,\ \ \rho'_2=(\beta'-1) \tan 2\beta.
\end{equation}
Then the curve $\gamma(t)$ whose polar coordinates are 
$$
\left(t-\beta(t)-\frac{\pi p}{q}, e^{\rho_1(t)}\right)
$$
is self-dual.
\end{proposition}

\paragraph{\it Proof.} Using the above  notation, the polar coordinates of the points $H(u(t))$ and $v(t)$ are 
$(t-\beta, e^{\rho_1})$ and $(t+\beta, e^{\rho_2})$. The differential equations (\ref{form}) are the Legendrian conditions $H(u)\cdot v'=0$ and $H(u')\cdot v=0$ which together imply that $H(u)\cdot v$ is constant. One needs to satisfy the monodromy conditions $\rho_1(t+c)=\rho_2(t), \rho_2(t+c)=\rho_1(t)$. Due to (\ref{form}), these equalities hold once one has
$$
0=\int_0^c \rho'_1(t)\ dt + \int_0^c \rho'_2(t)\ dt = 2\int_0^c \tan 2\beta\ d\beta=-\int_0^c d(\ln \cos 2\beta).
$$
The latter is zero because $\beta(c)=-\beta(0)$, and we are done. 
\proofend

 \paragraph{Radon curves.} Let $U$ be a Minkowski plane (two dimensional normed space) and let $\gamma$ be its unit circle, a closed smooth strictly convex centrally symmetric curve centered at the origin. For a vector $u\in \gamma$, one defines its orthogonal complement as the tangent line to $\gamma$ at $u$. This orthogonality relation is not symmetric, in general. A Minkowski plane is called a Radon plane, and the curve $\gamma$ a {\it Radon curve}, if the orthogonality relation is symmetric. A Radon curve admits a one-parameter family of circumscribed parallelograms whose sides are orthogonal to each other. Introduced by J. Radon about 90 years ago, Radon curves abound (they have functional parameters), see \cite{M-S} for a survey.
 
The relevance of Radon curves to our subject is the following statement. 

\begin{proposition} \label{Radon}
Radon curves are projectively self-dual.
\end{proposition}

\paragraph{\it Proof.} Let $[\ ,\ ]$ be an area element (linear symplectic structure) in $U$. Identify $U^\ast$ with $U$ using this area form:  $u^\ast = [\,\cdot ,u]$. With this identification,  $u^\perp$ is a line $\ell$, parallel to $u$, and such that $[v,u]=1$ for every $v\in \ell$. 

Similarly to the preceding discussion, a Radon curve can be realized as a curve $(u(t),v(t))$ in the space of pairs of non-zero vectors $(u,v)$, tangent to the distribution given by the conditions $[u,v']=0=[u',v]$ (these conditions mean that $u$ is orthogonal to $v$ and $v$ is orthogonal to $u$), and satisfying the monodromy conditions: $\displaystyle u\left(t+\frac\pi2\right)=v(t), v\left(t+\frac\pi2\right)=-u(t)$.  
The equalities $[u,v']=0=[u',v]$ imply that $[u,v]$ is constant or, after rescaling $\gamma$, that $[u,v]=1$. Therefore $u^\perp$ is the tangent line to $\gamma$ at point $v$, that is, $\gamma$ is self-dual.
\proofend

Let us conclude with two remarks. First, Radon planes can be also characterized as the Minkowski planes for which the unit circle is the solution to the isoperimetric problem (Busemann's theorem). Secondly, the outer billiard around a Radon curve possesses a one-parameter family of 4-periodic trajectories, see \cite{G-T} for a study of such outer billiards in the context of sub-Riemannian geometry.

 \end{document}